\documentclass{article}

\def\E{{\cal E}} \def\P{{\cal P}} \def\F{{\cal F}} \def\A{{\cal A}}
 \def\C{{\cal C}}  \def\K{{\cal K}}  
  \def\al{{\alpha}}

 \def\G{{\cal G}} \def\B{{\cal B}}
\newcommand\ord[1]{\circ (#1)}

\def\Sumer{\mathop{\sum\nolimits_r}\limits}

 \def\b{{\beta}}  

\title{Default Logic in a Coherent Setting}

\author{{\bf Giulianella Coletti}\\ Dipartimento   Matematica e Informatica \\Universit\`a di
Perugia,
  06100 Perugia (Italy) \\
 {\bf Romano Scozzafava}\\ Dipartimento
Metodi e Modelli Matematici\\   Universit\`a La Sapienza,   00161 Roma (Italy)
\\  {\bf Barbara Vantaggi}\\Dipartimento   Metodi e Modelli
Matematici\\  Universit\`a La Sapienza,  00161 Roma (Italy)\\ }

\begin{document}

\maketitle

\begin{abstract}
In this talk -- based on the results of a forthcoming paper (Coletti, Scozzafava and
Vantaggi 2002), presented also by one of us at the Conference on ``Non Classical
Logic, Approximate Reasoning and Soft-Computing'' (Anacapri, Italy, 2001) -- we
discuss the problem of representing default rules by means of a suitable coherent
conditional probability, defined on a  family of conditional events. An event is
singled-out (in our approach) by a  {\it proposition}, that is a statement that can be
either {\it true}  or  {\it false}; a conditional event is consequently defined by
means of two propositions and is  a 3--valued entity, the third value being (in this
context) a conditional probability.
\end{abstract}

\section{INTRODUCTION}

The concept of conditional event (as dealt with  in this paper) plays a central role
for the probabilistic reasoning. We give up  (or better, in a sense, we generalize)
the idea of de Finetti of looking at a conditional event $E|H$, with $H\neq\emptyset$
(the {\it impossible} event), as a $3$--valued logical entity looked on as
``undetermined'' when $H$ is false: it is {\it true} when both $E$ and $H$ are true,
{\it false} when $H$ is true and $E$ is false, while we let the third value {\it
suitably depend on the given ordered pair} $(E,H)$ and not being just an undetermined
{\it common value} for all pairs. It turns out (as explained in detail in Coletti and
Scozzafava 1999) that this function can be seen as a measure of the degree of belief
in the conditional event $E|H$, which under ``natural'' conditions reduces to the
conditional probability $P(E|H)$, in its most general sense related to the concept of
{\it coherence}, and satisfying the classic axioms as given by de Finetti (1949),
R\'enyi (1956), Krauss (1968), Dubins (1975): see Section 2. Notice that our concept
of conditional event differs from that adopted, for example, by Adams (1975),
Benferhat, Dubois and Prade (1997), Goodman and Nguyen (1988), Schay (1968).

Among the peculiarities (which entail a large flexibility in the management of any
kind of uncertainty) of this concept of {\it coherent} conditional probability versus
the usual one, we recall the following ones:\vspace{-.2cm}
\begin{itemize}
\item due to its {\it direct\/} assignment as a whole, the knowledge
(or the assessment) of the ``joint'' and ``marginal'' unconditional probabilities $P(E
\wedge H)$ and $P(H)$ is not required;\vspace{-.2cm}
\item the {\it conditioning} event $H$ (which {\it must} be a {\it possible} one) may
have {\it zero probability}, but in the assignment of $P(E|H)$ we are driven  by {\it
coherence}, contrary to what is done in those treatments where the relevant
conditional probability is given an {\it arbitrary} value in the case of a
conditioning event of zero probability;\vspace{-.2cm}
\item a suitable interpretation of its extreme values
$0$ and $1$ for situations which are different, respectively, from the trivial ones
$E\wedge H =\emptyset$ and $H\subseteq E$, leads to a ``natural'' treatment of the
{\it default reasoning}.
\end{itemize}\vspace{-.3cm}

In this talk we deal with the latter aspect.

\section{COHERENT CONDITIONAL \\ PROBABILITY}

The classic {\it axioms for a conditional probability} read as follows (given a set
$\C=\G \times \B^o$ of conditional events $E|H$  such that $\G$ is a Boolean algebra
and $\B \subseteq \G$ is closed with respect to (finite) logical sums, with
$\B^o=\B\setminus\{\emptyset\}\,$):\vspace{-.2cm}

\begin{itemize}
\item[{\it (i)\/}] $P(H|H) = 1$, for every $H\in \B^o\,$,

\item[{\it (ii)\/}] $P(\cdot|H)$ is a (finitely additive) probability on $\G$
for any given $H \in \B^o\,$,

\item[{\it (iii)\/}] $P(E \wedge A|H)=P(E|H)P(A|E \wedge H)$,\\
for any $A, E\in \G$,  $H, E \wedge H \in \B^o$.
\end{itemize}\vspace{-.2cm}

Conditional probability $P$ has been defined  on $\G\times\B^o$; however it is
possible, through the concept of {\it coherence}, to handle also those situations
where we need to assess $P$ on an {\it arbitrary} set ${\cal C}$ of conditional
events.

{\bf Definition 1} - The assessment $P(\cdot|\cdot)$ on ${\cal C}$ is {\it coherent}
if there exists $\C'\supset \C$, with $\C'= \G \times \B^o$, such that
$P(\cdot|\cdot)$ can be extended from $\C$ to $\C'$  as a {\it conditional
probability}.

A characterization of coherence is given (see, e.g., Coletti and Scozzafava 1996) by
the following

{\bf Theorem 1} - Let $\C$ be an arbitrary finite family of conditional events
$E_1|H_1, \ldots, E_n|H_n$ and $\A_o$ denote the set of  atoms $A_r$ generated by the
(unconditional) events $E_1,H_1,\ldots,E_n, H_n$. For a real function $P$ on $\C$ the
following two statements are equivalent:

(i) $P$ is a {\it coherent} conditional probability on $\C$;

(ii) there exists (at least) a {\it class} of probabilities $\{ P_0 ,P_1
,\ldots\,P_k\}$, each probability $P_\alpha$   being defined on a suitable  subset
$\A_\alpha \subseteq\A_o$, such that for any $E_i |H_i  \in \C$ there is a unique
$P_\alpha$ with
 $$ \Sumer_{A_r \subseteq H_i} P_\alpha (A_r) >0 \,,\,$$
 $$P(E_i|H_i) = \frac{\Sumer_{A_r \subseteq E_i  \wedge H_i} P_\alpha
(A_r)}{\Sumer_{A_r \subseteq H_i} P_\alpha (A_r)}\,\,;\leqno(1) $$ moreover
$\A_{\alpha '} \subset \A_{\alpha''}$ for $\alpha'
>\alpha''$ and $P_{\alpha''} (A_r) = 0$   if $A_r \in \A_{\alpha'}$.\vspace{1mm}

According to Theorem 1, a coherent conditional probability gives rise to a suitable
class $\{ P_o ,P_1 ,\ldots\,P_k\}$ of ``unconditional'' probabilities.

Where do the above classes of probabilities come from? Since $P$ is coherent on $\C$,
there exists an extension $P^*$ on $\G \times \B^o$, where $\G$ is the algebra
generated by the set $\A_o$ of atoms and $\B$ the additive class generated by $H_1,
\ldots, H_n$: then, putting $\F=\{\Omega, \emptyset\}$, the restriction of $P^*$ to
$\A_o \times \F^o$ satisfies (1) with $\alpha = 0$ for any $E_i|H_i$ such that
$P_o(H_i)>0\,$. The subset $\A_1\subset\A_o$ contains only the atoms $A_r \subseteq
H_o^1$, the union of $H_i$'s with $P_o(H_i)=0$ (and so on): we proved (see, e.g.,
Coletti and Scozzafava 1996, 1999) that, starting from a coherent assessment
$P(E_i|H_i)$ on $\C$, a relevant family $\P = \{P_\al\}$ can be suitably defined that
allows a representation such as (1). Every value $P(E_i|H_i)$ constitutes a constraint
in the construction of the probabilities $P_\al$ $(\al = 0, 1,...)$; in fact, given
the set $\A_o$ of atoms generated by $E_1,...,E_n,H_1,...,H_n$, and its subsets
$\A_{\alpha}$ (such that $P_\b(A_r )=0$ for  any $\b < \al$, with $A_r \in
\A_{\alpha}$), each $P_\al$ must satisfy the following system $(S_\alpha)$ with
unknowns $P_\alpha (A_r ) \ge 0$, $A_r \in \A_{\alpha}$, $$(S_\alpha) \cases{
\displaystyle \Sumer_{A_r \subseteq E_i H_i} P_\alpha (A_r) = P (E_i|H_i) \Sumer_{A_r
\subseteq H_i} P_\alpha (A_r)\,,\vspace{1mm} \,\, \cr
\smallskip \big[{\rm if}\ P_{\alpha-1} (H_i) =0\big]\cr \noalign{\medskip}
\displaystyle \Sumer_{A_r\subseteq H_0^\alpha} P_\alpha (A_r) =1} $$

\noindent where  $P_{-1} (H_i) =0$  for all $H_i$'s, and $H_o^\alpha$  denotes, for
$\alpha \ge 0$, the  union  of the $H_i$'s such that $P_{\alpha-1}(H_i ) = 0$; so, in
particular, $H_o^o  = H_o = H_1 \vee \ldots \vee H_n\,.$

Any class $\{ P_\alpha\}$ singled-out by the condition $(ii)$ is said {\it to agree}
with the conditional probability $P$.  Notice that in general there are infinite
classes of probabilities $\{P_\alpha\}$\,; in particular we have {\it only one
agreeing class} in the case that $\C$ is a product of Boolean algebras.

A coherent assessment $P$, defined on a  set $\C$ of conditional events, can be
extended in a natural way to all the conditional events $E|H$ such that $E \wedge H$
is an element of the algebra $\G$ spanned by the (unconditional) events $E_i, H_i\,,\,
i=1,2,...,n$ taken from the elements of $\C$, and $H$ is an element of the additive
class spanned by the $H_i$'s. Obviously, this extension is not unique, since there is
no uniqueness in the choice of the class $\{ P_\alpha\}$ related to condition {\it
(ii)} of  Theorem 1.

In general, we have the following result (see, e.g., Coletti and Scozzafava 1996):

{\bf Theorem 2} - If $\C$ is a  given  family of conditional events and  $P$  a
corresponding assessment,  then  there  exists  a (possibly not unique) coherent
extension of $P$ to an arbitrary family $\K$ of conditional events, with $\K \supseteq
\C$, {\it if and only if} $P$ is   coherent on $\C$.

Notice that if $P$ is coherent on a family $\C$, it is coherent also on $\E \subseteq
\C$.

\section{ZERO-LAYERS}

Given a class $\P = \{P_\alpha\}_{\alpha=0,1,\dots,k}$, agreeing with a conditional
probability on $\C$, it {\it naturally induces\/} the {\it zero-layer\/} $\ord{H}$ of
an event $H$, defined as \vspace{-2mm} $$\ord{H}=\beta\;\;\mbox{ if }
P_\beta(H)>0\,;\vspace{-2mm}$$ if $P_\alpha(H)=0$ for every $\alpha=0,1,\dots,k$
(obviously, we  necessarily have $H \neq H_i$ for every $i=1,2,\dots,n$\/), then
$\,\ord{H}=k+1$.

The zero-layer of a conditional event $E|H$ is defined as\vspace{-1mm}
$$\ord{E|H}=\ord{E\wedge H}-\ord{H}.\vspace{-1mm}$$ Obviously, for the certain event
$\Omega$ and for any event $E$ with positive probability, we have
$\ord{\Omega}=\ord{E}=0$ (so that, if the class contains only an {\it everywhere
positive} probability $P_o$, there is only one (trivial) zero-layer, {\it i.e.}
$\al=0$), while we put $\ord{\emptyset}=+\infty $. Clearly, \vspace{-1mm} $$\ord{A\vee
B}=\min\{\ord{A}, \ord{B}\}. \vspace{-.1cm}$$ Moreover, notice that $P(E|H)>0$ if and
only if $\ord{E H}= \ord{H}$, {\it i.e.} $\ord{E|H}=0$.

On the other hand, Spohn (see, for example, Spohn 1994, 1999) considers degrees of
plausibility defined via a {\it ranking} function, that is a map $\kappa$ that assigns
to each {\it possible} proposition a natural number (its {\it rank}) such
that\vspace{-.3cm}

\begin{enumerate}
\item[(a)] either $\kappa(A)=0$ or $\kappa(A^c)=0$, or both\,;
\item[(b)] $\kappa(A\vee B)=\min\{\kappa(A), \kappa(B)\}$\,;
\item[(c)] for all $A\wedge B\neq \emptyset$, the conditional rank of $B$ given
$A$ is \, $\kappa(B|A)=\kappa(A\wedge B)-\kappa(A)$\,.
\end{enumerate}\vspace{-.3cm}

Ranks represent degrees of ``disbelief''. For example, $A$ is {\it not} disbelieved
iff $\kappa(A)=0$, and it is disbelieved iff $\,\kappa(A)>0$.

{\bf Remark 1} - {\it Ranking functions are seen by Spohn as a tool to manage {\em
plain belief} and {\em belief revision}, since he maintains that probability is
inadequate for this purpose. In our framework this claim can be challenged (see {\em
Coletti, Scozzafava and Vantaggi 2001}), since our tools for belief revision are {\em
coherent conditional probabilities} and the {\bf ensuing} concept of {\em
zero-layers\/}: it is easy to check that zero-layers have the same formal properties
of ranking functions}.

\section{COHERENT PROBABILITY  AND DEFAULT LOGIC}

We recall that in Coletti, Scozzafava and Vantaggi (2001) we  showed that a sensible
use of events whose probability is $0$ (or $1$) can be a more general tool in revising
beliefs when new information comes to the fore, so that we have been able to challenge
the claim contained in Shenoy (1991) that probability is inadequate for revising plain
belief. Moreover, as recalled in Section 1, we may deal with the extreme value
$P(E|H)=1$ also for situations which are different from the trivial one $H\subseteq
E$.

The aim of  this Section is to handle, by means of a {\it coherent} conditional
probability, some aspects of {\it default reasoning} (see, e.g., Reiter 1980, Russel
and Norvig 1995): as it is well-known, a default rule is a sort of weak implication.

First of all, we discuss briefly some aspects of the classic example of Tweety.

The usual {\it logical implication} (denoted by $\subseteq$) can be anyway useful to
express that a penguin ($\pi$) is {\it certainly} a  bird ($\beta$),  i.e. $$\pi
\subseteq \beta\,,$$ so that $$P(\beta|\pi)=1\,;$$ moreover we know that Tweety
($\tau$) is a penguin (that is, $\tau \subseteq \pi$), and so also this fact can be
represented by a conditional probability equal to $1$, that is $$P(\pi|\tau)=1\,.$$

But we can express as well the statement ``a penguin {\it usually} does not fly'' (we
denote by $\varphi^c$ the contrary of $\varphi$, the latter symbol denoting
``flying'') by writing $$P(\varphi^c|\pi)=1\,.$$

(For simplicity, we have avoided to write down explicit a {\it proposition} -- that
is, an event -- such as ``a given animal is a penguin'', using the short-cut
``penguin'' and the symbol $\pi$ to denote this event; similar considerations apply to
$\beta$, $\tau$ and $\varphi$).

The question ``can Tweety fly?'' can be faced through an assessment of the conditional
probability $P(\varphi|\tau)$, which must be coherent with the already assessed ones:
by Theorem 1, it can be shown that {\it any value} $p\in [0,1]$ is a coherent value
for $P(\varphi|\tau)$, so that no conclusion can be reached -- {\it from the given
premises} -- on Tweety's ability of flying.

In other words, interpreting an equality such as $P(E|H)=1$ like a
{\it default rule} (denoted by $\longmapsto$), which in particular
(when $H \subseteq E$) reduces to the usual implication, we have
shown its {\it nontransitivity\,}: in fact we have
$$\tau\longmapsto\pi \mbox{\,\, and \,\,} \pi\longmapsto
\varphi^c\,,$$ but it {\it does not} necessarily follow the
further default rule $\tau\longmapsto \varphi^c$ (even if we {\it
might} have that $P(\varphi^c|\tau)=1$, {\it i.e.} that ``Tweety
usually does not fly'').

{\bf Definition 2} - Given a {\em coherent} conditional probability $P$ on a family
$\C$ of conditional events, a {\em default rule}, denoted by $H \longmapsto E$, is any
conditional event $E|H \in \C$ such that $P(E|H)=1$.

Clearly, any logical implication $A\subseteq B$ (and so also any equality $A=B$)
between events can be seen as a (trivial) default rule.

{\bf Remark 2} - {\it By resorting to the systems $(S_\alpha)$ to check the coherence
of the assessment $P(E|H)=1$ (which implies, for the relevant zero-layer,
$\ord{E|H}=0$), a simple computation gives $P_o(E^c \wedge H)=0$ (notice that the
class $\{P_\alpha\}$ has in this case only one element $P_o$). It follows
$\,\ord{E^c|H}=1$, so that $$\ord{E^c|H} > \ord{E|H}\,.$$  In terms of Spohn's ranking
functions (we recall -- {\em and underline} -- that our zero-layers are -- so to say
-- ``incorporated'' into a coherent conditional probability, so that {\bf we do not
need} an ``autonomous'' definition of ranking\,!) we could say, when $P(E|H)=1$, that
the disbelief in $E^c|H$ is greater than that in $E|H$. This conclusion {\bf must not}
be read as $P(E|H)> P(E^c|H)$\,!}\vspace{1mm}

Given a {\it set $\Delta \subseteq \C$ of default rules\/} $H_i \longmapsto E_i\,$,
with $i=1,...,n\,,$ we need to check its {\it consistency}, that is the coherence of
the ``global'' assessment $P$ on $\C$ such that $P(E_i|H_i)=1\,$,  $i=1,...,n\,$.

We stress that, even if our definition involves a conditional probability, the
condition given in the following theorem refers {\it only} to logical (in the sense of
{\it Boolean} logic) relations.

{\bf Theorem 3} - Given a {\em coherent} conditional probability $P$ on a family $\C$
of conditional events, the following two statements are equivalent:

(i) the set $\Delta \subseteq \C$ of default rules $$H_i \longmapsto E_i\,,\,\,
i=1,2,...,n\,,$$ represented by the assessment $$P(E_i|H_i)=1\,,\,\, i=1,2,...,n\,,$$
is consistent;

(ii) for every subset $$\{H_{i_1} \longmapsto E_{i_1},\ldots, H_{i_s} \longmapsto
E_{i_s} \}\,$$ of $\Delta\,$, with $s=1,2,...,n$, we have $$ \bigvee_{k=1}^s (E_{i_k}
\wedge H_{i_k}) \not\subseteq \bigvee_{k=1}^s (E_{i_k}^c \wedge H_{i_k}). \leqno(2)$$

{\it Proof} -  We prove that, assuming the above logical relations (2), coherence of
$P$ is compatible with the assessment $P(E_i|H_i)=1\;(i=1,2,...,n)\,,$ on $\Delta$.

We resort to the characterization Theorem 1: to begin with, put
$P(E_i|H_i)=1\;(i=1,2,...,n)\,,$ in the system $(S_o)$. The unconditional probability
$P_o$ can be obtained by putting $P_o(A_r)=0$ for all atoms \linebreak $A_r \subseteq
\bigvee_{j=1}^n\Big(E_j^c \wedge H_j\Big)$, so  for any atom $A_k \subseteq E_i \wedge
H_i $ which is not contained in $\bigvee_{j=1}^n\Big(E_j^c \wedge H_j\Big)\,$ --
notice that condition (2) ensures that there is such an atom $A_k$, since
$\bigvee_{j=1}^n (E_j\wedge H_j)\not\subseteq \bigvee_{j=1}^n (E_{j}^c \wedge H_{j})$
-- we may put $P_o(A_k)> 0$ in such a way that these numbers sum up to 1,  and we put
$P_o(A_r)=0$ for all remaining atoms.

This clearly gives a solution of the first system $(S_o)$. If, for some $i$, $E_i
\wedge H_i \subseteq \bigvee_{j=1 }^n\Big(E_j^c \wedge H_j\Big)$, then $P_o(E_i\wedge
H_i)=0$. So we consider the second system (which refers to all $H_i$ such that
$P_o(H_i)=0$), proceeding as above to construct the probability $P_1$; and so on.
Condition (2) ensures that at each step we can give positive probability $P_\alpha$ to
(at least) one of the remaining atoms.

Conversely, consider  the (coherent) assignment $P(E_i|H_i)=1$  (for $i=1,...,n$).
Then, for any index $j\in \{1,2,\dots,n\}$ there exists a probability $P_\alpha$ such
that $P_\alpha(E_{j} \wedge H_{j})>0$ and $P_\alpha(E_{j}^c \wedge H_j)=0$. Notice
that the restriction of $P$ to some conditional events
$E_{i_1}|H_{i_1},...,E_{i_s}|H_{i_s}$ of $\Delta$ is coherent as well.

Let $P_o$ be the first element of an agreeing class, and  $i_k$ an index such that
$P_o(H_{i_k})>0\,$: then we have  $P_o(E_{i_k} \wedge H_{i_k})>0$ and $P_o(E_{i_k}^c
\wedge H_{i_k})=0$. Suppose that $E_{i_k} \wedge H_{i_k} \subseteq \bigvee_{k=1}^s
(E_{i_k}^c \wedge H_{i_k})\,$: then $P_o(E_i{_k}\wedge H_{i_k})=0$. This contradiction
shows that condition (2) holds.\vspace{2mm}

{\bf Definition 3}  - A set $\Delta$ of default rules {\em entails\/} the default rule
$H\longmapsto E$ if the {\em only\/} coherent value for $P(E|H)$ is 1. In other words,
the rule $H\longmapsto E$ is entailed by $\Delta$ (or by a subset of $\Delta$) if
every possible extension (cf. Theorem 2) of the probability assessment
$P(E_{i_r}|H_{i_r})=1\,$, $r=1 \dots s\,,$ assigns the value $1$ also to $P(E|H)$.

Going back to the previous example of Tweety, its possible ability
(or inability) of flying can be expressed by saying that the
default rule $\tau \longmapsto \varphi$ (or $\tau \longmapsto
\varphi^c$)  {\it is not entailed} by the premises (the given set
$\Delta$).

\section{INFERENCE}

Several formalisms  for default logic have been studied  in the relevant literature
with the aim of discussing the minimal conditions that an entailment should satisfy.
In our framework this ``inferential'' process is ruled by the following

{\bf Theorem 4} - Given a set $\Delta$ of consistent default rules, we have:

{\bf (Reflexivity)}\vspace{-1mm}

\noindent $\Delta$ \,\, entails \,\, $A\longmapsto A\;\; \mbox{ for any }
A\neq\emptyset\,$;\vspace{2mm}

{\bf (Left Logical Equivalence)}\vspace{-1mm}

\noindent $(A=B) \,,\, (A\longmapsto C)\in \Delta$\, \,\, entails \,\, $ B\longmapsto
C\,$;\vspace{2mm}

{\bf (Right Weakening)}\vspace{-1mm}

\noindent $(A\subseteq B) \,,\, (C\longmapsto A)\in \Delta$\, \,\, entails \,\, $
C\longmapsto B\,$;\vspace{2mm}

{\bf (Cut)}\vspace{-1mm}

\noindent $(A\wedge B\longmapsto C) \,,\, (A\longmapsto B)\in \Delta$ \,\, entails
\,\, $A\longmapsto C\,$; \vspace{2mm}

{\bf (Cautious Monotonicity)}\vspace{-1mm}

\noindent $(A\longmapsto B) \,,\, (A\longmapsto C)\in \Delta$ \,\, entails \,\,
$A\wedge B\longmapsto C\,$; \vspace{2mm}

{\bf (Equivalence)}\vspace{-1mm}

\noindent $(A\longmapsto B) \,,\, (B\longmapsto A) \,,\, (A\longmapsto C)\in \Delta$
\,\, entails \,\, $ B\longmapsto C\,$; \vspace{2mm}

{\bf (And)}\vspace{-1mm}

\noindent $(A\longmapsto B) \,,\, (A\longmapsto C)\in \Delta$ \,\, entails \,\, $
A\longmapsto B\wedge C\,$; \vspace{2mm}

{\bf (Or)}\vspace{-1mm}

\noindent $(A\longmapsto C) \,,\, (B\longmapsto C)\in \Delta$ \,\, entails \,\, $
A\vee B\longmapsto C\,$.\vspace{2mm}

{\it Proof} - {\it Reflexivity\/} amounts to  $P(A|A)=1$ for every possible event.

 {\it Left Logical Equivalence\/} and {\it Right weakening\/} trivially
follow from elementary properties of conditional probability.

{\it Cut}\,: from $P(C|A\wedge B)=P(B|A)=1$  it follows that $$P(C|A)=P(C|A\wedge
B)P(B|A)+P(C|A\wedge B^c)P(B^c|A)=$$ $$=P(C|A\wedge B)P(B|A) =1\,.$$

{\it Cautious Monotonicity}\,: since $P(B|A)=P(C|A)=1$, we have that $$1= P(C|A\wedge
B)P(B|A)+
 P(C|A\wedge B^c)P(B^c|A)=$$ $$=P(C|A\wedge B)P(B|A)\,,$$
 hence $P(C|A\wedge B)=1$.

 {\it Equivalence}\,: since at least one conditioning event must have
positive probability, it follows that $A, B, C$ have positive probability; moreover,
$$P(A\wedge C)=P(A)=P(A\wedge B)=P(B)\,,$$ which implies $P(A\wedge B\wedge
C)=P(A)=P(B)$, so $P(C|B)=1$.

 {\it And}\,: since $$1\geq P(B\vee C|A)=P(B|A)+P(C|A)-P(B\wedge
C|A)=$$ $$=2-P(B\wedge C|A)\,,$$ it follows $P(B\wedge C|A)=1$.

 {\it Or}\,: since $$P(C|A\vee B)=$$
$$=P(C|A)P(A|A\vee B)+P(C|B)P(B|A\vee B)-$$ $$-P(C|A\wedge B )P(A\wedge B|A\vee B)=$$
$$= P(A|A\vee B)+P(B|A\vee B)-$$ $$-P(C|A\wedge B )P(A\wedge B|A\vee B)\geq 1\,,$$ we
get $P(C|A\vee B)=1$. \vspace{1mm}

We consider  now  some ``unpleasant'' properties (cf., e.g., Lehmann and Magidor,
1992), that in fact do not necessarily hold also in our framework:

{\bf (Monotonicity)}\vspace{-1mm}

\noindent $(A\subseteq B) \,,\, (B\longmapsto C) \in \Delta\,$ \,\, entails \,\, $
A\longmapsto C$ \vspace{2mm}

{\bf (Transitivity)}\vspace{-1mm}

\noindent $(A\longmapsto  B) \,,\, (B\longmapsto C) \in \Delta\,$ \,\, entails \,\, $
A\longmapsto C$ \vspace{2mm}

{\bf (Contraposition)}\vspace{-1mm}

\noindent $(A\longmapsto B) \in \Delta\,$ \,\, entails \,\, $ B^c\longmapsto
A^c$\vspace{2mm}

The previous example about Tweety  shows that {\it Transitivity\/} can fail.

In the same example, if we add the evaluation $P(\varphi|\beta)=1$ (that is, a bird
{\it usually\/} flies) to the initial ones, the assessment is still coherent (even if
$P(\varphi|{\pi})=0$ and $\pi\subseteq \beta$), but {\it Monotonicity\/} can fail.

Now, consider the conditional probability $P$ defined as follows\,: $$P(B|A)=1 \;,\;
P(A^c|B^c)=\frac{1}{4}\;;$$ it is easy to check that it is coherent, and so {\it
Contraposition\/} can fail.\vspace{1mm}

Many authors (cf., e.g., again Lehmann and Magidor, 1992) claim (and we agree) that
the previous unpleasant properties should be replaced by others, that we express below
in our own notation and interpretation: we show  that these properties hold in our
framework. Since a widespread consensus among their ``right'' formulation is lacking,
we will denote them as cs--(Negation Rationality),  cs--(Disjunctive Rationality),
cs--(Rational Monotonicity), where ``cs'' stands for ``in a coherent setting''. Notice
that, given a default rule $H\longmapsto E$, to say $(H\longmapsto E) \not\in \Delta$
means that the conditional event $E|H$ belongs to the set $\C
\setminus\Delta$.\vspace{2mm}

{\bf cs--(Negation Rationality)}

If $(A\wedge C \longmapsto B) \,,\, (A\wedge C^c \longmapsto B) \not\in \Delta$\\ then
$\Delta$ does not entail $ (A\longmapsto B)$\vspace{1mm}

{\it Proof} - If  $(A\wedge C \longmapsto B)$ and $(A\wedge C^c \longmapsto B)$ do not
belong to $\Delta$, i.e. $P(B|A\wedge C)<1$ and $P(B|A\wedge C^c)<1$  imply
$$P(B|A)=P(B|A\wedge C)P(C|A)+P(B|A\wedge C^c)P(C^c|A)<$$ $$<P(C|A)+P(C^c|A)=1\,.$$
\vspace{2mm}

{\bf cs--(Disjunctive Rationality)}

If $(A \longmapsto C) \,,\, (B \longmapsto C) \not\in\Delta$\\ then $\Delta$ does not
entail $(A\vee B\longmapsto C)$

{\it Proof} - Starting from the equalities $$P(C|A\vee B)=$$ $$= P(C|A)P(A|A\vee B)
+P(C|A^c\wedge B)P(A^c\wedge B|A\vee B)$$ and
 $$P(C|A\vee B)=$$
$$=P(C|B)P(B|A\vee B) +P(C|A\wedge B^c)P(A\wedge B^c|A\vee B),$$ since we have
$P(C|A)<1$ and $P(C|B)<1$, then  $P(C|A\vee B)=1$  would imply (by the first equality)
$P(A|A\vee B)=0$  and (by the second one) $P(B|A\vee B)=0$ (contradiction).
\vspace{2mm}

{\bf cs--(Rational Monotonicity)}

If $(A \wedge B\longmapsto C) \,,\, (A \longmapsto B^c) \not\in\Delta$\\ then $\Delta$
does not entail $ (A\longmapsto C)$

{\it Proof} - If it were $P(C|A)=1$, i.e. $$1=P(C|A\wedge B)P(B|A)+P(C|A\wedge
B^c)P(B^c|A)\,,$$  we would get either $$P(C|A\wedge B)=P(C|A\wedge B^c)=1$$ or one of
the following $$P(C|A\wedge B)=P(B|A)= 1 \,,$$ $$ P(C|A\wedge B^c)=P(B^c|A)= 1$$
(contradiction). \vspace{1mm}

In conclusion, let us notice the simplicity of our approach (Occam's razor...!), with
respect to other well-known methodologies, such as, e.g. those given  by Adams (1975),
Benferhat, Dubois and Prade (1997), Goldszmidt and Pearl (1996), Lehmann and Magidor
(1992), Schaub (1998).

\section{DISCUSSION}

Thought-provoking comments of two anonymous reviewers suggested to us to add this
further section.

Among coherence--based approaches to default reasoning (in the framework of imprecise
probability propagation), that of Gilio (2000) deserves to be mentioned, even if we
claim (besides the utmost simplicity of our definitions and results) many important
semantic and syntactic differences.

First of all, our framework (see the very beginning of our Introduction) is clearly
and rigorously settled: conditional events $E|H$ are {\bf not}  3-valued entities
whose third value is looked on as ``undetermined'' when $H$ is false, but they have
been defined instead in a way which entails ``automatically'' (so-to-say) the axioms
of conditional probability, which are those {\it ruling coherence\/}  (the details, as
already recalled in the Introduction, are in Coletti and Scozzafava, 1999).

In other words ({\bf french} words, since we are in France), ``tout se tient'', while
in the aforementioned paper by Gilio  a concept such as $E|H$ is interpreted sometimes
as a 3-valued entity looked on as ``undetermined'' when $H$ is false, sometimes as an
ordered pair of events, sometimes as a conditional assertion $H\,|\hspace{-2mm}\sim E$
(in the knowledge base).

Moreover, our notions of {\it consistency\/} and {\it
entailment\/} are both different from his: in fact he gives a
theorem (without proof) connecting the notion of consistency to
that of Adams (1975).

The problem is that we do not understand Adams' framework: in fact he requires
probability to be {\it proper\/} (i.e., positive) on the given events, but (since the
domain of a probability $P$ is an algebra)  we need to extend $P$ from the given
events to other events (by the way, coherence is nothing but complying with this
need). In particular, these ``new'' events may have zero probability: it follows,
according to Adams' definition of {\it conditional\/} probability in the case of a
conditioning event of zero probability, that we can easily get {\it incoherent\/}
assessments (see the example below). By the way, in the section ``{\bf Some
preliminaries}'', Gilio claims ``We can frame our approach to the problem of
propagating imprecise conditional probability assessments from the {\it probabilistic
logic\/} point of view, see, e.g., Frisch and Haddawy ...'': unfortunately, Frisch and
Haddawy definition of conditional probability coincides (for conditioning events which
are null) with that of Adams, and so it violates coherence as well!

Not to mention that both Gilio and Adams (and many others: some of them are mentioned
at the end of the previous section) base the concept of consistency on that of {\it
quasi conjunction\/}, which is a particular conditional event (and our concept of
conditional event is different from theirs); moreover we deem that the notion they
give of {\it verifiability\/} of a conditional event $E|H$, that is $E \wedge H \neq
\emptyset$, is too weak -- except in the case $H=\Omega$ -- to express properly the
relevant semantics.

Our discussion can be better illustrated by the following (very simple) example:

{\bf Example} - Consider two (logically independent) events $H_1$ and $H_2$, and put
$$E_1=H_1 \wedge H_2 \;,\; E_2=H_1^c \wedge H_2\,,$$ $$E_3=H_1^c \wedge H_2^c \;,\;
E=E_2 \;,\; H=H_3=\Omega\,.$$ Given $\alpha\,$, with $\,0<\alpha<1\,$, the assessment
$$P(E_1|H_1)=P(E_2|H_2)=1\;,\;P(E_3|H_3)=\alpha$$ on $\C=\{E_1|H_1 , E_2|H_2 ,
E_3|H_3\}$ is coherent; the relevant probabilities of the atoms are $$P(H_1 \wedge
H_2) = P(H_1 \wedge H_2^c)=0\,,$$ $$P(H_1^c \wedge H_2^c)= \alpha \,,\, P(H_1^c \wedge
H_2)=1-\alpha\,,$$ so that the set $\Delta$  of default rules corresponding to
$\{E_1|H_1 , E_2|H_2\}\,$ is consistent.

Does $\Delta$ entail $E|H\,$? A simple check shows that the only coherent assessment
for this conditional event is  $P(E|H)=1-\alpha$. Then the answer is NO, since we
require (in the definition of entailment) that $1$ is (the only) coherent extension.

On the contrary, according to Gilio characterization of entailment -- that is:
$\Delta$ (our notation) entails $E|H\,$ iff $P(E^c|H)=1$ is not coherent -- the answer
to the previous question is YES, since the only coherent value of this conditional
probability is $P(E^c|H)=\alpha\,$ (see the above computation).

For any $\epsilon>0$, consider now the assessment $$P(E_1|H_1)=1 \,,\,
P(E_2|H_2)=1-\epsilon\;,$$ so that $\{E_1|H_1 , E_2|H_2\}\,$ is consistent according
to Adams, as can be easily checked giving the atoms the probabilities $$P(H_1 \wedge
H_2) = \epsilon \,,\, P(H_1 \wedge H_2^c)=0\,,$$ $$P(H_1^c \wedge H_2^c)=0 \,,\,
P(H_1^c \wedge H_2)=1-\epsilon\,,$$  (notice that the assessment is proper). But for
any event $A \subset H_1 \wedge H_2^c$ we can extend $P$, according to his definition
of conditional probability, as $$P(A|H_1 \wedge H_2^c)=P(A^c|H_1 \wedge H_2^c)=1\,,$$
which is {\bf not} coherent!

Finally, there is no mention in Gilio's paper of Negation Rationality, Disjunctive
Rationality, and Rational Monotonicity (and, according to one of the reviewers, these
properties do not hold ``in default reasoning under coherent probabilities'', while in
our setting they have been proved at the end of Section 5).

\subsubsection*{Acknowledgements}
We thank an anonymous referee for signaling us a slight mistake in the proof of
Disjunctive Rationality.

\subsubsection*{References}

\quad\, E. Adams (1975). {\it The Logic of Conditionals}, Dordrecht: Reidel.

S. Benferhat, D. Dubois, and H. Prade (1997). Nonmonotonic Reasoning, Conditional
Objects and Possibility Theory. {\it Artificial Intelligence} {\bf 92}:259--276

G. Coletti and R. Scozzafava (1996). Characterization of Coherent Conditional
Probabilities as a Tool for their Assessment and Extension. {\it International Journal
of Uncertainty, Fuzziness and Knowledge-Based System} {\bf 4}:103--127.

G. Coletti and R. Scozzafava (1999). Conditioning and Inference in Intelligent
Systems. {\it Soft Computing} {\bf 3}:118--130.

G. Coletti, R. Scozzafava, and B. Vantaggi (2001). Probabilistic Reasoning as a
General Unifying Tool. In S.Benferhat and P. Besnard (eds.),  {\it Lectures Notes in
Computer Science} LNAI 2143 (ECSQUARU 2001), 120--131.

G. Coletti, R. Scozzafava, and B. Vantaggi (2002). Coherent Conditional Probability as
a Tool for Default Reasoning. {\it Proc. IPMU 2002}, Annecy (France), to appear.

B. de Finetti (1949). Sull'impostazione assiomatica  del  calcolo  delle
probabilit\`a. {\it Annali Univ. Trieste} {\bf 19}:3--55. (Engl. transl.:  Ch.5 In
{\it Probability, Induction, Statistics}, London: Wiley, 1972).

L.E. Dubins (1975). Finitely Additive Conditional Probabilities, Conglomerability and
Disintegration. {\it Annals of  Probability} {\bf 3}:89--99.

A. Gilio (2000). Precise propagation of upper and lower probability bounds in system
P. In {\it Proc. 8th Int. Workshop on Non-monotonic Reasoning}, ``Uncertainty
Frameworks in Non-Monotonic Reasoning'', Breckenridge (USA).

M. Goldszmidt and J. Pearl (1996) Qualitative probability for default reasoning,
belief revision and causal modeling. {\it Artificial Intelligence} {\bf 84}:57--112.

I.R. Goodman and H.T. Nguyen (1988). Conditional objects and the modeling  of
uncertainties. In M. Gupta and  T. Yamakawa  (eds.),  {\it Fuzzy Computing}, 119-138.
Amsterdam: North Holland.

P.H. Krauss  (1968). Representation of Conditional Probability Measures on Boolean
Algebras. {\it Acta Math. Acad. Scient. Hungar} {\bf 19}:229--241.

D. Lehmann and M. Magidor (1992).  What does a conditional knowledge base entail? {\it
Artificial Intelligence} {\bf 55}:1--60.

R. Reiter  (1980). A Logic for Default Reasoning. {\it Artificial Intelligence} {\bf
13}(1-2):81--132.

A. R\'enyi (1956). On Conditional Probability Spaces Generated by a Dimensionally
Ordered Set  of Measures. {\it Theory of Probability and its Applications} {\bf
1}:61--71.

S.J. Russel and P. Norvig (1995).  {\it Artificial Intelligence. A Modern Approach.}
New Jersey: Prentice-Hall.

T. Schaub (1998). The Family of Default Logics. In D.M. Gabbay and P. Smets (eds.),
{\it Handbook of Defeasible Reasoning and Uncertainty Management Systems}, Vol.2,
77--133, Dordrecht: Kluwer.

G. Schay (1968).  An Algebra of Conditional Events. {\it Journal of Mathematical
Analysis and Applications} {\bf 24}:334-344.

P.P. Shenoy (1991). On Spohn's Rule for Revision of Beliefs. {\it International
Journal of Approximate Reasoning} {\bf 5}:149--181.

W. Spohn (1994). On the Properties of Conditional Independence. In P. Humphreys, P.
Suppes (eds),  {\it Scientific Philosopher} {\bf 1}, Probability and Probabilistic
Causality, 173--194. Dordrecht: Kluwer.

W. Spohn  (1999). Ranking Functions, AGM Style. {\it Research Group
 ``Logic in Philosophy''},  Preprint 28.

\end{document}